\documentclass[prl,preprint]{revtex4}

\usepackage{graphicx,amssymb}
\usepackage{graphicx}
\usepackage{rotating}
\usepackage{dcolumn}
\usepackage{amsfonts}
\usepackage{amssymb}
\usepackage{amsmath}
\usepackage{latexsym}
\usepackage{epsfig}
\usepackage{dsfont,eucal,mathrsfs}
\usepackage{bm}

\usepackage{float}

\usepackage{color}

\renewcommand{\mathbf}{\boldsymbol}
\renewcommand{\mathcal}{\mathscr}

\definecolor{violet}{RGB}{111,0,255}

\begin{document}

\title{Improving the Network Structure can lead to Functional Failures}
\author{Jan Philipp Pade\footnote{pade@mathematik.hu-berlin.de}}
\affiliation{Humboldt University of Berlin, Institute of Mathematics, Unter den Linden 6, 10099 Berlin, Germany.}

\author{Tiago Pereira\footnote{tiago.pereira@imperial.ac.uk}}
    \affiliation{Department of Mathematics, Imperial College London, London SW72AZ, United Kingdom}
    \affiliation{London Mathematical Laboratory, London WC2N 6DF, UK}

\maketitle

{\bf 
In many real-world networks the ability to synchronize is a key property for its performance. Examples include  power-grid \cite{Motter2013},  sensor \cite{Sensor}, and neuron networks \cite{bookEp} as well as consensus formation \cite{Atay2013}. Recent work on undirected networks with diffusive interaction revealed that improvements in the network connectivity such as making the network more connected and homogeneous enhances synchronization \cite{Huang2008,Jalili2013,Motter2005}. However, real-world networks have directed and weighted connections \cite{Newman}.  In such directed networks, understanding the impact of structural  changes on the network performance remains a major challenge. Here, we show that improving the structure of a directed network can lead to a failure in the network function. For instance, introducing new links to reduce the minimum distance between nodes can lead to instabilities in the  synchronized motion. 
This counter-intuitive effect only occurs in directed networks. 
Our results allow to identify the dynamical importance of a link and thereby have a major impact on the design and control of directed networks.
}

Our everyday life depends on network synchronization at various levels. In power grids, power stations must keep a proper synchronization to avoid energy supply disturbances and blackouts \cite{Motter2013}. Sensor networks rely on synchronization among sensors to transmit information \cite{Sensor}. In the brain, epileptic seizures and Parkinson's diseases are a strong manifestation of synchronization \cite{bookEp}.  Further examples can be found in consensus formation \cite{Atay2013}. 

These complex systems are modeled by networks with diffusive interaction, that is,  the interaction between any two coupled  elements depends on the difference of their states. So far, research efforts to understand the influence of connectivity on the dynamics have focused on undirected diffusive networks. For instance,  it is known that increasing  the homogeneity or the number of connections enhances synchronization, as the maximum distance between nodes is decreased \cite{Huang2008,Jalili2013,Motter2005}.

Networks found in nature are often \textit{directed} and \textit{weighted}. For example, electrical synapses in neuron networks have asymmetric conductance \cite{Kandel}, which makes the underlying network directed.  Recent work has provided sufficient conditions to guarantee the stability of synchronization in directed networks in terms of the network structure and nature of the interaction \cite{Belykh,Pereira2014,MotterT}. However, understanding the impact of structural modifications, such as changing weights and adding or deleting links, on synchronization remains an open problem.

In this letter, using synchronization as a paradigm for network function, we show that improving the network connectivity structure can lead to a functional failure. This counter-intuitive phenomenon has dynamical origins. Namely, in directed networks the structural improvement has a suppressing effect in the network spectrum leading to the onset of instabilities associated with the synchronous motion. Furthermore, we identify a class of links in directed networks for which increasing the weights enhances synchronization. Our results provide a way to understand the dynamical importance of links and to develop strategies to avoid functional failures when improving the network structure. 

We consider directed networks of identical elements with diffusive interaction. 
The theory we develop here is general and can include networks of non-identical elements with minor modifications \cite{Pereira2014}. 
The network dynamics is described by
\begin{equation}\label{md1}
  \dot{\bm{x}}_i = \bm{f}(\bm{x}_i) + \alpha \sum_{j=1}^n W_{ij} \bm{H}(\bm{x}_j  - \bm{x}_i)\,,
\end{equation}
where $i=1,2,...,n$, $\bm{f} : \mathbb{R}^m \rightarrow \mathbb{R}^m$ is smooth, $\alpha\geq0$ is the overall coupling strength, and the matrix $\bm{W}$ describes the network structure, i.e., $W_{ij}\geq 0$ measures the strength of interaction from node $j$ to node $i$.  We assume that the network solutions are bounded and  the local coupling function $\bm{H}$ is smooth satisfying $\bm{H}(\bm{0}) = \bm{0}$. This last condition guarantees that the synchronous state $\bm{x}_1 = \bm{x}_2 = \cdots =\bm{x}_n$ is a solution of the coupled equations for all values of $\alpha.$   The overall coupling strength $\alpha$ represents the fixed energy cost per connection.

\begin{figure*}[ht!] 
   \centering
   \includegraphics[width=4.5in]{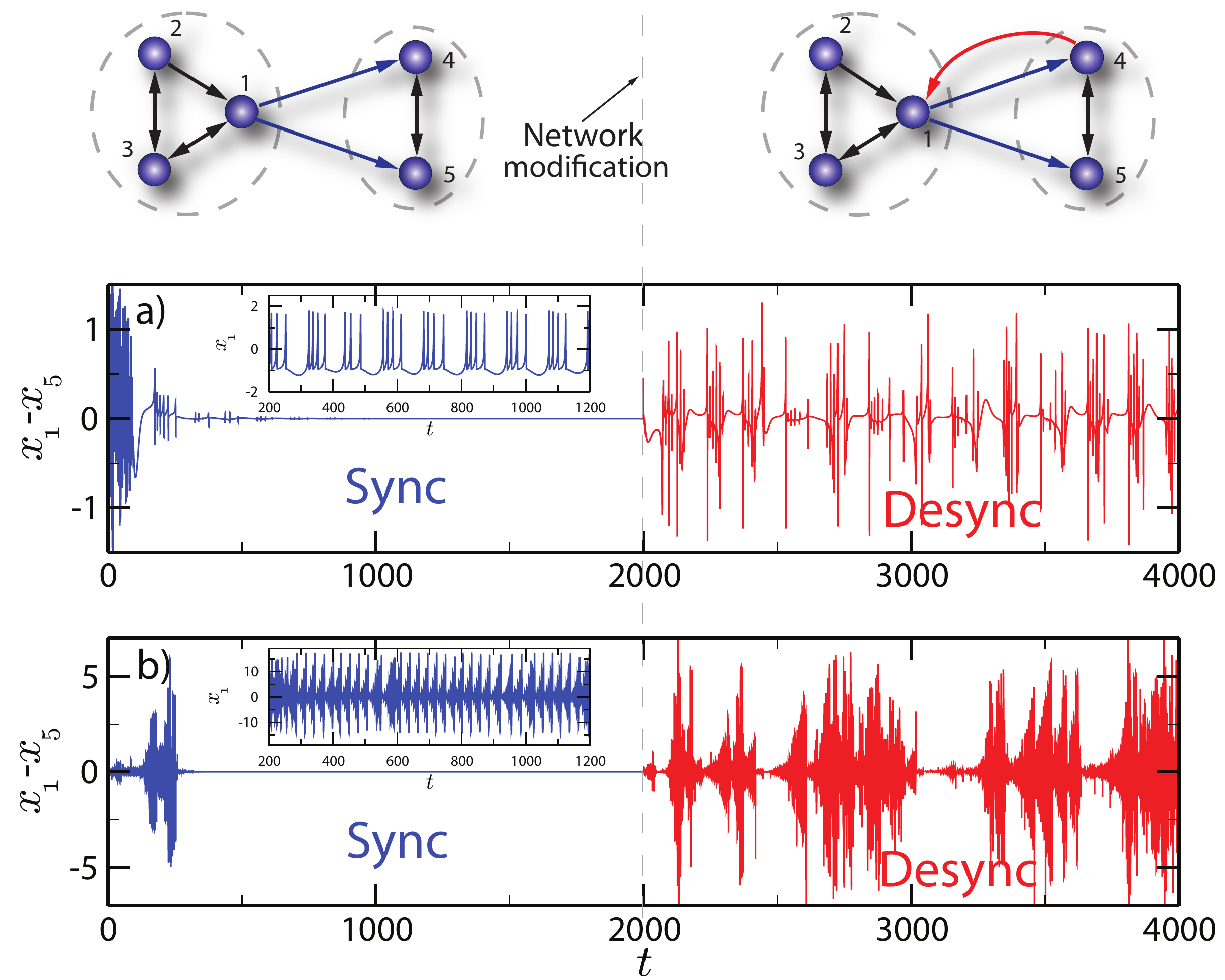} 
   \caption{Improving connectivity leads to desynchronization. The figures show simulation results for the networks on top. All links in blue have weight one. In the main plots we show the difference of the first component of $x_1$ and the first component of $x_5$. 
In a) the node dynamics is given by Hindmarsh-Rose (HR) neurons, and in b) by Roessler dynamics. The global coupling $\alpha$ is chosen such that the nodes synchronize chaotically for the original network. This can be seen in the main plots for times until $t=2000$ in blue. The introduction of the new link $4\rightarrow 1$ with weight $0.4$ at time $t=2000$ leads to a destabilization, displayed in red. The insets show the time series of a single node. For the HR neurons we consider a chaotic bursting mode and for the Roessler dynamics a chaotic state. }
   \label{fig:DesyncSimulation}
\end{figure*}

Here, we consider \textit{weakly connected} networks, that is, when ignoring the link's directions the network is connected.
A directed network is \textit{strongly connected} if every node can be reached by every other node through a directed path. If a directed network is not strongly connected the links can be partitioned into two different classes: links belonging  to some strongly connected subnetwork, called \textit{strongly connected component}, and links belonging to some \textit{cutset}. A cutset is a set of links which point from one strongly connected component to another \cite{Bang-Jensen2009}.
In the top left inset of Figure \ref{fig:DesyncSimulation} we show a network composed of two strongly connected subnetworks indicated by the grey dotted ellipses, and in blue the cutset.

The smaller strongly connected component does not influence the larger component, as there are no links from the smaller to the larger component. Nonetheless, the network still supports stable synchronous dynamics. The counter-intuitive effects of directed networks come into play once we try to improve the network structure. Introducing a new link pointing from node $4$ to node $1$ improves the connection structure significantly, as the whole network is now strongly connected: any two nodes in the network are connected by a directed path. However, this structural improvement has a surprising consequence for the dynamics: the synchronous state becomes unstable.

We illustrate this effect with two different classes of node dynamics, namely,  Hindmarsh-Rose neurons in the chaotic bursting regime \cite{Barrio2011} as shown in the inset of Figure \ref{fig:DesyncSimulation} a), and chaotic Roessler oscillators \cite{Roessler1979}, see inset of Figure \ref{fig:DesyncSimulation} b). The state of each node is given by a three-dimensional vector $\bm{x}_i = (x_i,y_i,z_i)$. Details on the models can be found in the supplementary material. All nonzero weights $W_{ij}$ in the network on top left are set to one and we choose the global coupling $\alpha$ such that the nodes synchronize, that is, for any nodes $i$ and $j$ the difference of states $\bm{x}_i(t) - \bm{x}_j(t)$ vanishes for $t\rightarrow\infty$. 
This synchronous dynamics can be seen in Figure \ref{fig:DesyncSimulation} for times $t<2000$.  At time $t=2000$ we add the new link $4\rightarrow 1$ with a weight of $W_{14}=0.4$, which leads to the strongly connected network on the top right. As can be seen for times $t>2000$ this destabilizes the synchronous state.

\begin{figure}[h!] 
   \centering
   \includegraphics[width=3in]{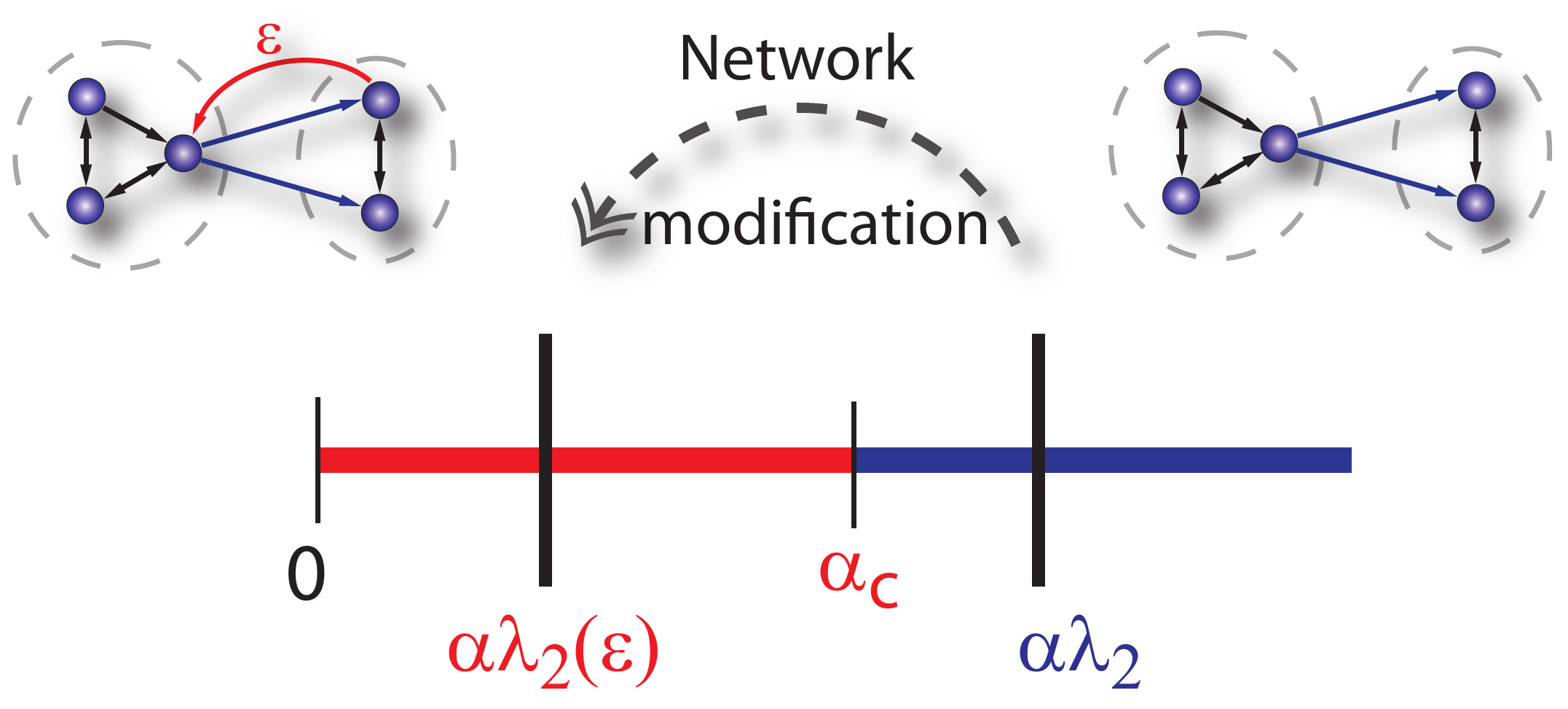} 
   \caption{Motion of the spectral gap. In a schematic representation we illustrate the motion of the spectral gap $\lambda_2$ under structural modifications in case $\lambda_2(\varepsilon)$ is real. The network on the right has a spectral gap such that $\alpha\lambda_2>\alpha_c$. Adding a link as indicated decreases the gap to $\lambda_2(\varepsilon)$, which violates the stability condition (\ref{StabCond}).}
   \label{fig:EigenvalueShift}
\end{figure}

To understand this phenomenon, we analyze the stability of the synchronization subspace $\bm{\gamma}(t)=\bm{x}_1(t)=\bm{x}_2(t)=\cdots=\bm{x}_n(t)$, with $\dot{\bm{\gamma}}(t)=\bm{f}(\bm{\gamma}(t))$. The variational equation of Eq. (\ref{md1}) along $\bm{\gamma}(t)$ can be decomposed into $n$ blocks of the form
\begin{equation}\label{decoupled_system}
\dot{\bm{\xi}}_i = \left[D\bm{f}(\bm{\gamma}(t)) -  \alpha \lambda_i\bm{\Gamma}\right]\bm{\xi}_i,
\end{equation}
where $D\bm{f}$ denotes the Jacobian,  $\bm{\Gamma} = D\bm{H}(\bm{0})$. The $\lambda_i$ are the eigenvalues  of the \textit{Laplacian matrix} $\bm{L}=\bm{D}-\bm{W}$ where $\bm{D}$ is the diagonal matrix with the row sums of $\bm{W}$ on its diagonal \cite{Pereira2014,MotterT}. We assume that the spectrum of $\bm{\Gamma}$ is real, which is the case for many applications. Now, these are decoupled $m$-dimensional equations which only differ by the Laplacian eigenvalues $\lambda_i$. We consider the case where $\bm{L}$ has a simple zero eigenvalue \cite{Ren2007}.  In this case, the eigenvalues of $\bm{L}$ have nonnegative real parts, so we can order them increasingly according to their real parts 
$0 = \lambda_1 < \Re(\lambda_2) \le \Re(\lambda_3) \le \cdots \le \Re(\lambda_n)$.  
The eigenvalues with non-zero real parts correspond to dynamics transverse to the synchronization subspace. Therefore, if the corresponding equations (\ref{decoupled_system}) have stable trivial solutions, synchronization in Eq. (\ref{md1}) is stable. For a large class of coupling functions and local dynamics \cite{Pereira2014}, the stability condition for synchronization is given by  
\begin{equation}\label{StabCond}
\alpha \Re(\lambda_2) > \alpha_c
\end{equation}
where $\alpha_c = \alpha_c(\bm{f},\bm{\Gamma})$ (see supplementary material for more details).  More involved stability conditions can be tackled, but the analysis becomes more technical without providing new insight into the phenomenon.
Condition (\ref{StabCond}) shows that the \textit{spectral gap} $\lambda_2$ plays a central role for synchronization properties of the network. Structural changes which decrease the real part of $ \lambda_2$ can destabilize the synchronous state, see Figure \ref{fig:EigenvalueShift}.

Let us consider the generic case where the eigenvalue $\lambda_2$ is simple. Using perturbation analysis \cite{Lancaster1985}, we obtain the direction of growth of $\lambda_2$ as a function of the structural modifications in the network. This analysis leads to a characterization of links which are capable of destabilizing the network. The spectral gap $\lambda_2(\varepsilon)$ of a perturbed Laplacian $\bm{L}_p =  \bm{L}+\varepsilon\tilde{\bm{L}}$ is given by  $\lambda_2(\varepsilon) = \lambda_2 + \lambda_2^{\prime} \varepsilon + O(\varepsilon^2) $ with  
\begin{equation}\label{Eigenv_Pert}
\lambda_2^{\prime} =  \frac{\langle \bm{u}, \tilde{\bm{L}} \bm{v} \rangle  }{\langle \bm{u}, \bm{v} \rangle}
\end{equation}
where $\langle \cdot, \cdot \rangle$ is the Euclidean inner product, $\bm{u}$  a left and $\bm{v}$ a right eigenvector of the unperturbed Laplacian $\bm{L}$ corresponding to the spectral gap $\lambda_2$, and we choose the vectors such that $\langle \bm{u} ,\bm{v}\rangle>0$.

If the network is undirected then $\bm{L}$ is symmetric and the left and right eigenvectors are dual.  As $\tilde{\bm{L}}$ is positive semi-definite we obtain $\lambda_2^{\prime} \geq 0$.  So increasing weights as well as adding links do not decrease the spectral gap in undirected networks \cite{Fiedler1973}. This is an essential difference to \textit{directed} networks. 

Let us now consider a directed network composed of two strongly connected components as in the example in Figure \ref{fig:DesyncSimulation}. Our approach is general and this choice is for the sake of simplicity. 
In this situation the Laplacian can be represented in block form
$$
\bm{L}=\left(\begin{array}{cc}
\bm{L}_{1} & \bm{0}\\
-\bm{C} & \bm{L}_{2}+\bm{D}_c
\end{array}\right)
$$ 
where $\bm{C}$ represents the cutset pointing from the strong component $\bm{W}_1$ to $\bm{W}_2$, $\bm{L}_{1,2}$ are the respective Laplacians and $\bm{D}_c$ is again a diagonal matrix with the row sums of $\bm{C}$ on its diagonal. As a consequence of the block structure, eigenvalues of $\bm{L}$ are either eigenvalues of $\bm{L}_1$ or eigenvalues of $\bm{L}_2+\bm{D}_c$. 
Suppose the eigenvalue with the smallest nonzero real part is located in the second component $\bm{W}_2$ (this encloses the example from Figure \ref{fig:DesyncSimulation}). Using a Perron-Frobenius argument,  we can show that the eigenvalue $\lambda_2$ is real, see supplementary material. So, the eigenvectors $\bm{u}$ and $\bm{v}$ are real and by equation (\ref{Eigenv_Pert}) the motion of $\lambda_2$ is along the real axis. 

To determine the effects of  the network changes on the spectral gap we investigate the left and right eigenvectors $\bm{u} = (\bm{p},\bm{q})$ and $\bm{v} = (\bm{r},\bm{s})$. This decomposition corresponds to the triangular form of the Laplacian.  Using the block structure, we can show that $\bm{r}=\bm{0}$ and that $\bm{q}$ and $\bm{s}$ are left and right eigenvectors of $\bm{L}_{2}+\bm{D}_c$. Furthermore, again by a Perron-Frobenius argument we can show that both $\bm{q}$ and $\bm{s}$ are positive.
By our assumptions, $\lambda_2$ is not an eigenvalue of $\bm{L}_{1}$. Therefore, we can solve the eigenvector equation for $\bm{p}$ in terms of $\bm{q}$ to obtain $\bm{p}=\bm{q}\bm{C}(\bm{L}_{1}-\lambda_2 \bm{I})^{-1}$, where $\bm{I}$ is the identity matrix.

Now we introduce a link in opposite direction of the cutset, so the Laplacian writes as
\[
\bm{L}_p =\left(\begin{array}{cc}
\bm{L}_{1}+\bm{D}_{\Delta} & -\bm{\Delta}\\
-\bm{C} & \bm{L}_{2}+\bm{D}_c
\end{array}\right),
\]
where $\bm{\Delta}$ is the matrix describing the new link, and $\bm{D}_{\Delta}$ is the associated diagonal matrix, see supplementary material. This yields
\begin{equation}\label{EigenvPert2}
\lambda_2^{\prime} = -\frac{\langle \bm{q}, \bm{M} \bm{s} \rangle}{\langle \bm{q}, \bm{s} \rangle},
\end{equation}
where 
$$
\bm{M} = \bm{C}(\bm{L}_1-\lambda_2 \bm{I} )^{-1}\bm{\Delta}
$$
takes into account the structural changes $\bm{\Delta}$. Now, such a modification will weaken the stability or even lead to instabilities if $\lambda_2^{\prime} < 0$. Determining the modifications that yield a decrease of the spectral gap is an involved problem, and we shall tackle it elsewhere. Here, we will focus on the example of Fig. 1, as it contains all central concepts without technical intricacies.  

In the example from Figure \ref{fig:DesyncSimulation}, $\bm{L_2}+\bm{D}_c$ is symmetric and we have $\bm{q}=\bm{s}=-\frac{1}{\sqrt{2}}(1,1)$. Moreover,  $\bm{M} =\frac{1}{2}\left(\begin{array}{cc}
0 & 1\\
0 & 1
\end{array}\right)$ and $\bm{q}$ and $\bm{s}$ are eigenvectors of $\bm{M}$ with eigenvalue $1/2$.
Because $\bm{q}$ is a common eigenvector of both $\bm{L_2}+\bm{D}_c$ and $\bm{M}$ corresponding to a positive eigenvalue we obtain a decrease of the spectral gap with a rate  $\lambda_2^{\prime}=-\frac{1}{2}$. 
This is a main mechanism that generates instabilities: {\it the eigenvectors of $\bm{L_2}+\bm{D}_c$ lie in the space spanned by the eigenvectors of $\bm{M}$ with positive eigenvalues.} 

Formally, we can obtain all the structural changes capable for destabilization as a function of the eigenvectors $\bm{q},\bm{s}$ of $\bm{L}_1$ and $\bm{L}_2+\bm{D}_c$. In contrast to undirected networks where there is a well developed theory relating eigenvectors to the underlying graph structure \cite{Biyikoglu2007}, for directed graphs the theory is underdeveloped.  Therefore, further analytical insights remain a challenge. From a computational point of view though, we can solve this problem for any given network. 

For modifications $\bm{\Delta}$ in the direction of the cutset the situation is clear, as Eq. (\ref{Eigenv_Pert}) reduces to $\lambda_2^{\prime}=\frac{\langle \bm{q} , \bm{D}_{\Delta}\bm{s} \rangle}{\langle \bm{q} , \bm{s} \rangle }$. Because all the involved quantities are positive, the spectral gap does not decrease when reinforcing the cutset: it represents a \textit{stabilising} class of links. By increasing strengths in the cutset we are guaranteed to enhance synchronization.

Our results revealed that directed and undirected networks behave essentially distinct under structural changes. Namely, assuming the stability condition Eq. (\ref{StabCond}), synchronization loss caused by structural improvements is a property inherent to directed networks exclusively. 
If $\alpha \Re(\lambda_2) \gg \alpha_c$, the network modification may not destroy synchronization. However, it worsens the quality in the sense that the transient towards synchronization becomes larger. 

Recently, interconnected networks have attracted much attention  \cite{Buldyrev2010,Radicchi2013}, as they can exhibit catastrophic cascades of failures when connections are undirected. Our results suggest that interconnected networks in which interconnections are represented by (directed) cutsets behave qualitatively different.

The catastrophic effects of structural improvements on the network function have a long history in game theory.
In the realm of games such effects are known as Braess's paradox \cite{Braess2005}. In games the effect occurs because the players take rational decisions to optimise their strategies. In the case of complex networks of dynamical systems considered here, the effect is dynamical and is a consequence of the motion of eigenvalues of the network Laplacian. Furthermore, our results shed light on how to plan and design network modifications without destroying the network performance, as for instance discussed for power-grids in \cite{Timme2012,Motter2013}.

{\bf Acknowledgments:} We are in debt with D. Turaev, R. Medrano, S. Yanchuk for valuable discussions. 
This work was partially supported by FAPESP-DFG International Research Training Group (IRTG) 1740 and
Marie Curie IIF Fellowship (303180).

\appendix

\section{Supplementary Material}

\section{Simulations}

All numerics were done with Matlab2011a using a fourth order Runge-Kutta method. We remind that the main equations are given by
\begin{equation}\label{main_equation}
\dot{\bm{x}}_i=\bm{f}(\bm{x}_i)+\alpha	\sum_{j=1}^na_{ij}\bm{H}(\bm{x}_j-\bm{x}_i).
\end{equation}

\subsection{Hindmarsh-Rose oscillators}

The Hindmarsh-Rose model is a three dimensional ordinary differential equation which models the membrane potential of a neuron. Depending on the parameter settings it exhibits spiking and bursting behaviour. For this model, the local dynamics $\bm{f}$ is given by
\begin{eqnarray*}
\dot{x}&=&y+a_1x^2-x^3-z+I \\
\dot{y}&=&1-5x^2-y\\
\dot{z}&=&a_2(s(x-x_R)-z).
\end{eqnarray*}
Here, $I$ is a constant input current. The parameters are chosen as follows: $a_1=3.01$, $a_2=0.006$, $s=4$, $I=3.2$ and $x_R=-1.6$. In this regime we can observe chaotic bursting in the local dynamics  \cite{Barrio2011}. We consider the electrical synaptic interaction between neurons given by 
\[
\bm{H}=\left(\begin{array}{ccc}
1 &0 &0 \\
0 &0 &0 \\
0 &0 &0 
\end{array}\right),
\]
so the local coupling is only in the $x$-component, known as membrane potential. 
The stability condition $\alpha\Re(\lambda_2)>\alpha_c$ can be verified via a master stability function approach \cite{Pecora1998}. 
In order to achieve stable synchronized motion for the whole network we fixed $\alpha = 0.96$.

\subsection{Roessler oscillators}

For the Roessler oscillators the local dynamics are given by
\begin{eqnarray*}
\dot{x}&=&-y-z \\
\dot{y}&=&x+a_1y\\
\dot{z}&=& a_2+z(x-a_3).
\end{eqnarray*}

where we chose $a_1=0.2$, $a_2=0.2$ and $a_3=9$. We consider the interaction in all variables 
\begin{eqnarray*}
\bm{H}=\left(\begin{array}{ccc}
1 &0 &0 \\
0 &1 &0 \\
0 &0 &1 
\end{array}\right).
\end{eqnarray*}
In this case, applying the results from \cite{Pereira2014} we again obtain the stability condition $\alpha\Re(\lambda_2)>\alpha_c$. In order to achieve stable synchronisation we fixed $\alpha =0.092$.

\section{Laplacians of digraphs and associated eigenvectors}
Suppose we have the following Laplacian
$$
\bm{L}=\left(\begin{array}{cc}
\bm{L}_1 & \bm{0}\\
-\bm{C} & \bm{L}_2+\bm{D}_{\bm{C}}
\end{array}\right),
$$ 
for two strongly connected components, where $\bm{C}\neq \bm{0}$ represents the cutset pointing from the strongly connected component associated to $\bm{L}_1\in\mathbb{R}^{n_1\times n_1}$ to the one associated to $\bm{L}_2\in\mathbb{R}^{n_2\times n_2}$. We assume that the underlying network has a rooted spanning tree, then the zero eigenvalue is simple \cite{Chebotarev2000}. Because of the block form of the Laplacian, the characteristic polynomial factorizes to
\[
p(\lambda)=\det(\lambda-\bm{L}_1)\det(\lambda-(\bm{L}_2+\bm{D}_{\bm{C}})).
\]
Now, suppose the nonzero eigenvalue with smallest real part $\lambda_2$ is simple and given as a zero of the second factor, so it is an eigenvalue of $\bm{L}_2+\bm{D}_{\bm{C}}$. Then we can show that this eigenvalue is real and positive and the corresponding left and right eigenvectors of $\bm{L}_2+\bm{D}_{\bm{C}}$ are positive as well. We begin by showing that the eigenvectors are positive. To do so, consider $s=\max_{i}\left\{ d_{i}+\sum_{j\neq i}a_{ij}\right\}$, then $\bm{N}=s\bm{I}-\left(\bm{L}_2+\bm{D}_{\bm{C}}\right)$ is a nonnegative matrix by definition of $s$. Furthermore, as the component associated to $\bm{L}$ is strongy connected, $\bm{N}$ is irreducible. Then, by the Perron-Frobenius theorem \cite{Lancaster1985}, $\bm{N}$ has a maximal real eigenvalue $\Lambda$ with corresponding nonnegative left and right eigenvectors $\bm{\omega}$ and $\bm{\eta}$. That is 
\begin{eqnarray*}
\bm{N}\bm{\eta} &=&\Lambda\bm{\eta}\\
\iff\left(\bm{L}_2+\bm{D}_{\bm{C}}\right)\bm{\eta}	&=&	\left(s-\Lambda\right)\bm{\eta} .
\end{eqnarray*}
As $\Lambda$ is the maximal eigenvalue and all the eigenvalues of $\bm{L}_2+\bm{D}_{\bm{C}}$ are obtained by eigenvalues $\mu$ of $\bm{N}$ through $s-\mu$, we must have that $s-\Lambda$ is real and the minimal eigenvalue of $\bm{L}_2+\bm{D}_{\bm{C}}$. Furthermore, $\bm{N}$ and $\bm{L}_2+\bm{D}_{\bm{C}}$ have the same eigenvectors, so the left and right eigenvectors corresponding to $s-\Lambda$ are nonnegative, which proves the first statement. For the positiveness of the minimal eigenvalue $s-\Lambda$ of $\bm{L}_2+\bm{D}_{\bm{C}}$ first remark that $\bm{L}_2+\bm{D}_{\bm{C}}$ is diagonally dominant
\begin{eqnarray*}
(\bm{L}_2+\bm{D}_{\bm{C}})_{ii}&=&(\bm{D}_{\bm{C}})_{ii}+\sum_{k\neq i}(\bm{L}_2)_{ik}\\
&\geq &\sum_{k\neq i}(\bm{L}_2)_{ik}\\
&=&\sum_{k\neq i}(\bm{L}_2+\bm{D}_{\bm{C}})_{ik}.
\end{eqnarray*}
By the Gershgorin theorem \cite{Lancaster1985} every eigenvalue lies in at least one of the circles with center $(\bm{L}_2+\bm{D}_{\bm{C}})_{ii}$ and radius $\sum_{k=1}^n(\bm{L}_2+\bm{D}_{\bm{C}})_{ik}$. Consequently, the eigenvalues cannot have negative real parts. As by our assumptions  the zero eigenvalue is simple, $s-\Lambda$ has to be positive.

\end{document}